\documentclass[12pt,a4paper,reqno]{amsart}

\usepackage{amsfonts,amsthm,amsmath,amssymb,amscd,mathrsfs}
\usepackage{indentfirst}
\usepackage{cite}
\usepackage{bm,xcolor}
\usepackage{enumitem}
\usepackage{mathtools}
\usepackage{graphicx}
\usepackage{tikz}
\usetikzlibrary{arrows.meta,decorations.markings}
\usepackage{caption}
\usepackage[hidelinks]{hyperref} 

\theoremstyle{plain} 
\newtheorem{theorem}{Theorem}[section]
\newtheorem{remark}{Remark}[section]
\newtheorem{lemma}[theorem]{Lemma}
\newtheorem{prop}[theorem]{Proposition}

\theoremstyle{definition}

\def\mbN{\mathbb{N}}
\def\mbR{\mathbb{R}}
\def\mbS{\mathbb{S}}
\def\mbT{\mathbb{T}}

\def\bv{\mathbf{v}}
\def\bw{\mathbf{w}}
\def\bn{\mathbf{n}}

\def\bT{\mathbf{T}}

\def\mcC{\mathcal{C}}

\def\mcL{\mathcal{L}}

\def\mcT{\mathcal{T}}

\def\msC{\mathscr{C}}

\def\msN{\mathscr{N}}
\def\msR{\mathscr{R}}

\def\Div{{\rm div}}

\newcommand{\Rn}[1]{\romannumeral #1\relax}

\newcommand{\btau}{\boldsymbol{\tau}}
\newcommand{\bga}{\boldsymbol{\gamma}}

\begin{document}

\title[A characterization of shear flows in an infinite strip]{A necessary and sufficient condition for steady Euler flows to be shear flows in an infinite strip}

\author{Changfeng Gui}
\address{Department of Mathematics, Faculty of Science, University of Macau, Taipa, Macao.}
\email{changfenggui@um.edu.mo}

\author{Xiaofan Wang}
\address{Department of Mathematics, Faculty of Science, University of Macau, Taipa, Macao.}
\email{yc47905@um.edu.mo}

\author{Huan Xu}
\address{Department of Mathematics, Faculty of Science, University of Macau, Taipa, Macao.}
\email{hzx0016@auburn.edu}

\author{Wen Yang}
\address{Department of Mathematics, Faculty of Science, University of Macau, Taipa, Macao.}
\email{wenyang@um.edu.mo}



\date{}

\begin{abstract}
We give a necessary and sufficient characterization of parallel shear flows among {\it all} $C^1$ steady incompressible Euler flows in an infinite strip, subject to the slip boundary condition. We show that a steady Euler flow is a parallel shear flow if and only if almost every regular streamline traverses the strip longitudinally. The proof relies on three key ingredients. First, the topological assumption on the streamlines yields global bounds on the stream function, the velocity, and the vorticity. Second, we establish a novel geometric formula for the total curvature which, together with the topological assumption on the streamlines, implies the finiteness of the total curvature via the area and coarea formulas. Third, once the total curvature is finite, we derive the asymptotic behavior of the flow in the far field, which in turn forces the total curvature to vanish.
\end{abstract}

\maketitle



\section{Introduction}

This paper studies the two-dimensional steady incompressible Euler equations
\begin{eqnarray}\label{Euler}
\left\{\begin{aligned}
&\bv\cdot\nabla\bv+\nabla P=0, & \hbox{in } \Omega, \\
&\Div\,\bv=0, & \hbox{in } \Omega,
\end{aligned}\right.
\end{eqnarray}
where $\bv=(v_1,v_2)$ is the velocity field, $P$ is the scalar pressure, and $\Omega$ is a planar domain with a smooth boundary $\partial\Omega$. 
The flow is subject to the slip boundary condition
\begin{align}\label{SBC}
\bv\cdot\bn=0\ \hbox{ on } \partial\Omega,
\end{align}
where $\bn$ denotes the outward unit normal to $\partial\Omega$.

Among steady solutions of the Euler system, two basic classes are parallel shear flows (or simply shear flows) and circular flows. By a parallel shear flow we mean a velocity field which, after a rotation of coordinates, has the form $\bv=(v_1(x_2),0)$; by a circular flow we mean a velocity field of the form $\bv=V(|x|)\frac{x^\perp}{|x|}$, where $x^\perp=(-x_2,x_1)$.

We are concerned with the rigidity of steady Euler flows. A basic question is to determine under what conditions a steady flow must be a parallel shear flow or a circular flow. The systematic study of this question owes much to the pioneering work of Hamel and Nadirashvili \cite{HN_CPAM_2017,HN_ARMA_2019,HN_JEMS_2023}. In \cite{HN_CPAM_2017}, they proved that a $C^2$ steady Euler flow in an infinite strip, satisfying the slip boundary condition and $\inf |\bv|>0$, must be a parallel shear flow. They subsequently established the corresponding rigidity theorem in the whole plane for bounded flows whose speed has a positive lower bound \cite{HN_ARMA_2019}. In annular domains, they \cite{HN_JEMS_2023} also established circular symmetry under suitable assumptions on the stagnation set. Their work also developed a close connection between Euler rigidity and symmetry questions for semilinear elliptic equations, which has been central to many subsequent developments (see, e.g., \cite{DER_CMP_2026,DN_CPAM_2026,EHSX_Duke_2026,EFR_arXiv_2025,FWXZ_arXiv_2026,GPSY_Duke_2021,LLSX_Poincare_2026,Ruiz_ARMA_2023,WZ_arXiv_2023}).

Our main interest here is the rigidity of shear flows. Most existing results give sufficient conditions for rigidity, such as a smallness condition \cite{CEW_ARMA_2023,LZ_ARMA_2011} or a stagnation-free condition \cite{HN_CPAM_2017,HN_ARMA_2019}. These conditions have led to strong rigidity theorems, but they need not be satisfied by every shear flow. It is therefore natural to seek conditions which are both necessary and sufficient, while allowing general stagnation sets.

Recent work of Gui, Xie and Xu \cite{GXX_CMP_2026,GXX_arXiv_2026} provides such characterizations. In the whole plane, they proved in \cite{GXX_CMP_2026} that a bounded smooth steady Euler flow is a parallel shear flow if and only if its velocity direction omits at least one point of $\mbS^1$. Thus, every bounded nonshear flow in the plane must attain every velocity direction, even when stagnation points are present. In a periodic channel $\mbT\times(-1,1)$, Gui, Xie and Xu \cite{GXX_arXiv_2026} obtained a stronger, purely topological characterization: a steady Euler flow is a shear flow if and only if it has no contractible closed streamlines.

In the present paper, we consider the infinite strip
\[
\Omega_\infty:=\mbR\times(-1,1).
\]
In this domain, the absence of closed streamlines alone does not imply shear rigidity. Indeed, there exist bounded nonshear steady Euler flows without stagnation points whose streamlines turn back towards the same end of the strip \cite{DR_CVPDE_2025,GRXX_ARMA_2026}. The relevant distinction is therefore whether a regular streamline joins the two different ends of the strip. This brings us back to the work of Hamel and Nadirashvili \cite{HN_CPAM_2017}. An essential part of their proof is to show that every streamline traverses the strip from one end to the other. Although their theorem is stated under a sufficient nondegeneracy condition, their analysis already identifies the geometric property which underlies shear rigidity in a strip. Our purpose is to prove that this property is itself necessary and sufficient. In fact, it suffices to impose it on almost every regular level of the stream function.

To state the result precisely, let $\bv\in C^1(\overline{\Omega_\infty})$ be a divergence-free velocity field satisfying the slip boundary condition \eqref{SBC}, and define the stream function by
\[
u(x_1,x_2)=\int_{x_2}^{1}v_1(x_1,t)\,dt.
\]
Then $u\in C^2(\overline{\Omega_\infty})$ satisfies
\[
\bv=\nabla^\perp u:=(-\partial_{x_2}u,\partial_{x_1}u).
\]
Moreover, $u(x_1,1)=0$, and the slip boundary condition gives
\[
u(x_1,-1)=c
\]
for some constant $c\in\mbR$. We denote by $\msC(u)$ the set of critical values of $u$, that is,
\[
\msC(u)=\{ u(x): x\in \overline{\Omega_\infty},\, |\nabla u(x)|=0 \}.
\]
Then $\msR(u)=u(\overline{\Omega_\infty})\setminus\msC(u)$ denotes the set of regular values of $u$, where $u(\overline{\Omega_\infty})$ is the image of $\overline{\Omega_\infty}$ under $u$. A connected component of $u^{-1}(r)$, with $r\in\msR(u)\setminus\{0,c\}$, is a regular streamline. Such a streamline is called {\it longitudinal} if it connects the two infinite ends of the strip. More precisely, it admits a parametrization $\bga=(\gamma_1,\gamma_2):\mbR\to\Omega_\infty$ such that, after reversing the parametrization if necessary,
\[
\gamma_1(s)\longrightarrow \pm \infty \quad \hbox{as } s\to\pm\infty.
\]
Let $\mcL^1$ denote the Lebesgue measure on $\mbR$.

Our main result is stated as follows.

\begin{theorem}\label{thm_main}
Let $\bv\in C^1(\overline{\Omega_\infty})$ and $P\in C^1(\overline{\Omega_\infty})$ solve \eqref{Euler}--\eqref{SBC}. Then the following statements are equivalent:
\begin{enumerate}[label=\textup{(\roman*)}]
\item $\bv$ is a shear flow, that is, $\bv$ is of the form $\bv=(v_1(x_2),0)$.

\item For every $r\in\msR(u)\setminus\{0,c\}$, every connected component of $u^{-1}(r)$ is longitudinal.

\item There is an $\mcL^1$-null set $\msN\subset\mbR$ such that, for every $r\in\msR(u)\setminus(\msN\cup\{0,c\})$, every connected component of $u^{-1}(r)$ is longitudinal.
\end{enumerate}  
\end{theorem}

A few remarks on the theorem are in order.

\begin{remark}
The implications (\Rn{1}) $\Rightarrow$ (\Rn{2}) $\Rightarrow$ (\Rn{3}) are evident. Under assumption (\Rn{2}) or (\Rn{3}), $\msC(u)$ may be nonempty. Consequently, the stagnation set may have positive area and may contain open subsets.
\end{remark}

\begin{remark}
The implication (\Rn{3}) $\Rightarrow$ (\Rn{1}) shows that, for any nonshear flow, $\msR(u)$ has a subset of positive measure such that for every $r$ in this subset, at least one component of $u^{-1}(r)$ fails to be longitudinal. Such a component is either an interior Jordan curve or a properly embedded line with both ends tending to the same horizontal infinity.
\end{remark}

\begin{remark}
The same topological condition does not imply rigidity for steady Euler flows in the half-plane. Indeed,
\[
\bv=(\cosh x_1,-x_2\sinh x_1),
\qquad
P=\frac12(x_2^2-\cosh^2x_1)
\]
solves \eqref{Euler}--\eqref{SBC} in $\mbR\times(0,\infty)$, and $v_1\ge1$. Its stream function is $u=-x_2\cosh x_1$, and its interior level sets $u=-r$, $r>0$, are the graphs of the bounded functions
\[
x_2=\frac{r}{\cosh x_1},\qquad r>0.
\]
All these streamlines connect the two horizontal ends, while the flow is not a shear flow.
\end{remark}

The proof of Theorem~\ref{thm_main} is based on total curvature estimates. This method was first introduced into the study of steady Euler equations in \cite{GXX_CMP_2026}, and was further developed in \cite{GXX_arXiv_2026,GLXX_arXiv_2026}. For a subdomain $D\subset\Omega$, the total curvature is
\begin{equation*}
\mcT(\bv,D)
:=\int_D\bigl(|\nabla\bv|^2-|\nabla|\bv||^2\bigr)\,dx
=\int_D\frac{|v_1\nabla v_2-v_2\nabla v_1|^2}{|\bv|^2}\,dx,
\end{equation*}
where the quotient is defined to be zero on $\{\bv=0\}$. Away from stagnation points, if $\theta$ is a local angle of the velocity, then the integrand equals $|\bv|^2|\nabla\theta|^2$. Thus the total curvature measures the variation of the velocity direction, and its vanishing provides a way to prove shear rigidity even in the presence of stagnation points.

Several ways to estimate this quantity have been developed. In \cite{GXX_CMP_2026,GXX_arXiv_2026}, elliptic estimates for the vorticity equation in the form
\begin{equation*}
\Div\,(v_1\nabla v_2-v_2\nabla v_1)=0
\end{equation*}
play a central role. Another approach in \cite{GXX_CMP_2026} uses the coarea formula for the direction field. More recently, Gui, Li, Xie and Xu \cite{GLXX_arXiv_2026} obtained a representation of the total curvature by applying the coarea formula to pressure level sets. The main new ingredient here is a different formula, which allows us to estimate the total curvature by integrating along the level sets of the stream function, namely, along the streamlines themselves. This formula connects the total curvature directly with the topological assumption in Theorem~\ref{thm_main}.

The starting point is an identity obtained in the periodic-channel argument of \cite{GXX_arXiv_2026}. Under the assumption that there are no contractible closed streamlines, the authors proved
\[
\int_{\mbT\times(-1,1)}\big(|\nabla\bv|^2-|\nabla|\bv||^2\big)\,dx
=-\frac12\int_{\mbT\times(-1,1)}\Delta P\,dx.
\]
The right-hand side vanishes by periodicity and the slip boundary condition. Motivated by this identity, we find the following formula
\begin{equation}\label{eq_intro-geometric-identity}
|\nabla\bv|^2-|\nabla|\bv||^2=-\frac12\Delta P+\omega\bv\cdot\nabla\theta.
\end{equation}
Here $\theta$ is the local angle of the velocity, so
\[
\nabla\theta=\frac{v_1\nabla v_2-v_2\nabla v_1}{|\bv|^2}.
\]

The usefulness of \eqref{eq_intro-geometric-identity} comes from the fact that $\omega$ is constant along each regular streamline. For a bounded substrip $D=(a,b)\times(-1,1)$, the coarea formula gives
\begin{equation}\label{eq_intro-streamline-coarea}
\int_D\omega\bv\cdot\nabla\theta\,dx
=\int_{\mbR}\sum_{\Gamma\in\mcC_r(D)}\omega_\Gamma\,\Delta_\Gamma\theta\,dr,
\end{equation}
where $\mcC_r(D)$ denotes the collection of connected components of $u^{-1}(r)\cap D$, $\omega_\Gamma$ is the constant vorticity on $\Gamma$, and $\Delta_\Gamma\theta$ is the signed change of the velocity angle along $\Gamma$, oriented by the flow. The cancellation in this signed angular increment is essential. Simple geometric bounds control $|\Delta_\Gamma\theta|$ for arcs cut out by the vertical sides of $\partial D$, and the area formula controls the number of these arcs through the velocity on the two cross sections.

The paper is organized as follows. In Section~\ref{sec_global-bounds}, we establish the global bounds for the stream function, velocity and vorticity that follow from the longitudinal behavior of almost every regular streamline. In Section~\ref{sec_finiteTC}, we prove the finiteness of the total curvature by means of \eqref{eq_intro-geometric-identity} and \eqref{eq_intro-streamline-coarea}. Section~\ref{sec_asymptotic} concerns the far-field behavior of the pressure and velocity. Finally, in Section~\ref{sec_vanishTC}, we prove that the total curvature vanishes and complete the proof of Theorem~\ref{thm_main}.

\medskip

\section{Global bounds}\label{sec_global-bounds}

In this section, we assume condition (\Rn{3}) of Theorem~\ref{thm_main} and prove the global boundedness of the stream function, the velocity, and the vorticity. Here, the vorticity 
\[
\omega:=\partial_{x_1}v_2-\partial_{x_2}v_1=\Delta u
\]
is continuous in $\overline{\Omega_\infty}$ and solves
\[
\Div(\omega\bv)=0
\]
in the sense of distributions. It is easy to see that $\omega$ is constant along every regular streamline. In fact, $\omega\bv$ also admits a stream function $h\in C^1$ such that $\omega\bv=\nabla^\perp h$. Then $\nabla h=\omega\nabla u$, since $\bv=\nabla^\perp u$. Near any non-stagnation point, using $u$ as a local coordinate gives $h=f(u)$ for some $f\in C^1$. Hence, $\omega=f'(u)$. It means that $\omega$ is locally constant on streamlines and therefore constant along every regular streamline.

\begin{lemma}\label{lemma_global-bound}
Under condition (\Rn{3}) of Theorem~\ref{thm_main}, it holds that
\begin{align*}
\|u\|_{L^\infty(\Omega_\infty)}+\|\bv\|_{L^\infty(\Omega_\infty)}+\|\omega\|_{L^\infty(\Omega_\infty)}<\infty.
\end{align*}
Moreover, for every $p\in(1,\infty)$, it holds that
\begin{align*}
\sup_{s\in\mbR}\|u\|_{W^{2,p}((s-1,s+1)\times(-1,1))}<\infty.
\end{align*}
\end{lemma}
\begin{proof}
First, we show that
\[
m_0:=\min_{[-1,1]}u(0,\cdot)\le u \le M_0:=\max u(0,\cdot) \quad \hbox{ in } \Omega_\infty.
\]
By the classical Morse-Sard theorem, the critical values of the $C^2$ stream function $u$ form an $\mcL^1$-null subset of $\mbR$. If $u(y)>M_0$ at some point $y\in\Omega_\infty$, continuity of $u$ along a path from $y$ to $x_1=0$ would yield a regular value $r\in\msR(u)\setminus(\msN\cup\{0,c\})$ with $r>M_0$. Every component of $u^{-1}(r)$ is longitudinal and therefore meets $x_1=0$, a contradiction. The lower bound is identical.

Next, we use a similar argument to show the vorticity bound
\[
|\omega|\le W_0:=\max|\omega(0,\cdot)| \quad \hbox{ in } \Omega_\infty.
\]
If $|\omega(z)|>W_0$ for some $z\in\Omega_\infty$, then by continuity we have for some $\delta>0$ that $|\omega|>W_0$ in $B_{\delta}(z):=\{x:|x-z|<\delta\}$. 
In particular, $|\omega|>0$ in $B_{\delta}(z)$, which implies that the stream function $u$ is nonconstant in $B_{\delta}(z)$, since otherwise $\omega=\Delta u=0$.
It follows that $u(B_{\delta}(z))$ is a nonempty interval, hence contains a regular value $r\in\msR(u)\setminus(\msN\cup\{0,c\})$.
Now $u^{-1}(r)$ has a connected component $\Gamma$ intersecting the disk $B_{\delta}(z)$. Since the vorticity is conserved on $\Gamma$, one has $|\omega|>W_0$ on $\Gamma$. On the other hand, $\Gamma$ also intersects the section $\{0\}\times(-1,1)$ by assumption, so $|\omega|\le W_0$ on $\Gamma$, which is a contradiction. 
This proves the global pointwise vorticity bound.

Introduce the function $\Tilde{u}=u+\frac12 c(x_2-1)$ that solves
\[
\Delta\Tilde{u}=\omega, \quad \Tilde{u}(\cdot,\pm1)=0.
\]
The standard elliptic estimates give
\begin{align*}
\|\Tilde{u}\|_{W^{2,p}((s-1,s+1)\times(-1,1))}
\le C_p\big(\|\omega\|_{L^p((s-2,s+2)\times(-1,1))}+\|\Tilde{u}\|_{L^p((s-2,s+2)\times(-1,1))}\big),    
\end{align*}
with the constant $C_p$ depending only on $p$.
Since both $\omega$ and $\Tilde{u}$ are bounded in $\Omega_\infty$, the right side is uniformly bounded in $s\in\mbR$. 
The local $W^{2,p}$-estimates for $u$ then follow readily.
Finally, taking $p>2$ and using Sobolev embedding bounds $\nabla u$, hence $\bv$. This completes the proof of the lemma.
\end{proof}


\medskip

\section{Finiteness of the total curvature}\label{sec_finiteTC}

In this section, under assumption (\Rn{3}) of Theorem~\ref{thm_main} and relying on the global boundedness established previously, we show that the total curvature is finite.

\begin{prop}\label{prop_finiteTC}
Let $\bv\in C^1(\overline{\Omega_\infty})$ and $P\in C^1(\overline{\Omega_\infty})$ solve \eqref{Euler}--\eqref{SBC} such that assumption (\Rn{3}) of Theorem~\ref{thm_main} is fulfilled. Then it holds that
\begin{align}\label{ineq_TCbound}
\mcT(\bv,\Omega_\infty)
\le 4\pi\|\bv\|_{L^\infty(\Omega_\infty)}\|\omega\|_{L^\infty(\Omega_\infty)}.
\end{align}
\end{prop}

The key ingredient is the following pointwise formula, inspired by a global formula in \cite{GXX_arXiv_2026}.

\begin{lemma}\label{lemma_key-identity}
Let $\bv\in C^1(\overline{\Omega_\infty})$ and $P\in C^1(\overline{\Omega_\infty})$ solve \eqref{Euler}--\eqref{SBC}. It holds that
\begin{align}\label{eq_key-identity}
\frac{|v_1\nabla v_2-v_2\nabla v_1|^2}{|\bv|^2}=-\frac12\Delta P+\omega\bv\cdot\frac{v_1\nabla v_2-v_2\nabla v_1}{|\bv|^2}, \quad \hbox{ a.e. in } \Omega_\infty.
\end{align}    
\end{lemma}
\begin{proof}
The first and last terms in \eqref{eq_key-identity} are understood to vanish on the stagnation set, so both belong to $L^{\infty}_{\rm loc}(\overline{\Omega_\infty})$. However, the regularity of $P$ requires clarification. Applying the divergence operator to the momentum equation in \eqref{Euler} and using the incompressibility condition, we obtain
\[
-\Delta P=\partial_{x_1}\bv\cdot\nabla v_1+\partial_{x_2}\bv\cdot\nabla v_2=2((\partial_{x_1}v_1)^2+\partial_{x_1}v_2\partial_{x_2}v_1)\in C(\overline{\Omega_\infty}).
\]
Consequently, by elliptic interior regularity, $P\in W^{2,p}_{\rm loc}(\Omega_\infty)$ for every $p\in[1,\infty)$. In fact, owing to the slip boundary condition, $(\bv,P)$ can be extended to a locally Lipschitz solution of \eqref{Euler} in the whole plane by an elementary reflection (see \cite{GLXX_arXiv_2026}). Hence the local boundary regularity reduces to interior regularity, and $P\in W^{2,p}_{\rm loc}(\overline{\Omega_\infty})$.

On the other hand, a straightforward computation yields
\begin{align*}
-((\partial_{x_1}v_1)^2+\partial_{x_1}v_2\partial_{x_2}v_1)
=&\partial_{x_1}v_1\partial_{x_2}v_2-\partial_{x_1}v_2\partial_{x_2}v_1\\
=&\nabla^\perp v_1\cdot\nabla v_2\\
=&\frac{v_1^2\nabla^\perp v_1\cdot\nabla v_2-v_2^2\nabla^\perp v_2\cdot\nabla v_1}{|\bv|^2}\\
=&(v_1\nabla^\perp v_1+v_2\nabla^\perp v_2)\cdot\frac{v_1\nabla v_2-v_2\nabla v_1}{|\bv|^2}.
\end{align*}
Consequently,
\[
\frac12\Delta P=(v_1\nabla^\perp v_1+v_2\nabla^\perp v_2)\cdot\frac{v_1\nabla v_2-v_2\nabla v_1}{|\bv|^2}.
\]
We proceed as follows
\begin{align*}
-\frac12\Delta P+\omega\bv\cdot\frac{v_1\nabla v_2-v_2\nabla v_1}{|\bv|^2}
=&\frac{v_1\nabla v_2-v_2\nabla v_1}{|\bv|^2}\cdot(\omega\bv-v_1\nabla^\perp v_1-v_2\nabla^\perp v_2)\\
=&\frac{|v_1\nabla v_2-v_2\nabla v_1|^2}{|\bv|^2}.
\end{align*}
This gives \eqref{eq_key-identity} and completes the proof of the lemma.
\end{proof}

For a regular simple $C^1$ arc $\Gamma$ parametrized by $\bga=\bga(t):[0,1]\to\mbR^2$, its unit tangent vector 
\[
\btau=\frac{\bga'}{|\bga'|}
\]
admits a continuous lifting $\theta:[0,1]\to\mbR$ such that
\[
\btau(t)=(\cos\theta(t),\sin\theta(t)).
\]
Any two continuous liftings differ by a constant in $2\pi\mathbb Z$, so up to sign, the endpoint difference
\[
\Delta_{\Gamma}\theta:=\theta(1)-\theta(0)
\]
is independent of the choice of continuous liftings.

To estimate the second term on the right-hand side of \eqref{eq_key-identity}, we need the following elementary geometric lemma.

\begin{lemma}\label{lemma_angle-bound}
Let $\Gamma$ be a regular simple $C^1$ arc parametrized by $\bga=\bga(t):[0,1]\to\mbR^2$, and let $\theta\in C([0,1])$ be a continuous lifting of $\btau:=\bga'/|\bga'|$.
\begin{enumerate}[label=\textup{(\roman*)}]
\item
If a straight line $\ell$ contains $\bga(0),\bga(1)$ and $\Gamma$ lies on one side of $\ell$, then
\begin{equation*}
|\Delta_{\Gamma}\theta|\le 2\pi.
\end{equation*}
More precisely, in coordinates in which $\bga(0)=(0,a)$, $\bga(1)=(0,b)$, $a<b$, and $\Gamma\subset\{x_1\ge0\}$, let $\alpha\in[-\frac{\pi}{2},\frac{\pi}{2}]$ and $\beta\in[\frac{\pi}{2},\frac{3\pi}{2}]$ be the respective angle representatives of $\btau(0)$ and $\btau(1)$. Then
\begin{equation*}
\theta(1)-\theta(0)=\beta-\alpha\in[0,2\pi].
\end{equation*}

\item
If $\ell_0,\ell_1$ are distinct parallel straight lines, $\bga(0)\in\ell_0$, $\bga(1)\in\ell_1$, and $\Gamma$ lies between these lines, then
\begin{equation*}
|\Delta_{\Gamma}\theta|\le\pi.
\end{equation*}
More precisely, in coordinates in which the strip is $\{0\le x_1\le h\}$, $h>0$, with $\bga(0)=(0,a)$ and $\bga(1)=(h,b)$, let $\alpha,\beta\in[-\frac{\pi}{2},\frac{\pi}{2}]$ be the respective angle representatives of $\btau(0)$ and $\btau(1)$. Then
\begin{equation*}
\theta(1)-\theta(0)=\beta-\alpha\in[-\pi,\pi].
\end{equation*}
\end{enumerate}   
\end{lemma}
\begin{proof}
Let
\[
D=\{(s,t)\in\mbR^2:0\le s\le t\le1\}.
\]
For $(s,t)\in D$, define
\begin{equation*}
\bT(s,t)
=\begin{cases}
\frac{\bga(t)-\bga(s)}{|\bga(t)-\bga(s)|},
&s<t,\\[2mm]
\btau(s),&s=t.
\end{cases}
\end{equation*}
Then $\bT$ is a continuous map from $D$ to $\mbS^1$. Since $D$ is simply connected, $\bT$ admits a continuous lifting $\Theta$. In particular, $\Theta(t,t)$ is also a continuous lifting of $\btau$. Hence,
\[
\theta(1)-\theta(0)=\Theta(1,1)-\Theta(0,0).
\]

(\Rn{1}) 
Use the coordinates in part~(\Rn{1}). Then $\Theta(0,1)$ is the argument of the vector $\bT(0,1)=(0,1).$ We may choose $\Theta(0,1)=\pi/2$, which uniquely determines $\Theta$. For $0\le t\le 1$, the vector $\boldsymbol{\gamma}(t)-\boldsymbol{\gamma}(0)$ has nonnegative horizontal component, so $\mathbf{T}(0,t)$ lies in the closed right semicircle. By continuity of $\Theta(0,t)$ and the specified value $\Theta(0,1)=\pi/2$, we have $\Theta(0,t)\in[-\pi/2,\pi/2]$; in particular, $\Theta(0,0)\in[-\pi/2,\pi/2]$. Likewise, for $0\le s\le 1$, the vector $\boldsymbol{\gamma}(1)-\boldsymbol{\gamma}(s)$ has nonpositive horizontal component, so $\mathbf{T}(s,1)$ lies in the closed left semicircle. By continuity of $\Theta(s,1)$ and the specified value $\Theta(0,1)=\pi/2$, we have $\Theta(s,1)\in[\pi/2,3\pi/2]$; in particular, $\Theta(1,1)\in[\pi/2,3\pi/2]$. This proves part~(\Rn{1}).

(\Rn{2}) 
Use the coordinates in part~(\Rn{2}). Both $\bga(t)-\bga(0)$ and $\bga(1)-\bga(s)$ have nonnegative horizontal components. Thus, both unit vector fields $\bT(0,t)$ and $\bT(s,1)$ take values in the closed right semicircle. Now $\Theta(0,1)$ is the argument of the vector $(h,b-a)$ with $h>0$. We may choose $\Theta(0,1)\in(-\pi/2,\pi/2)$. Then both $\Theta(0,t)$ and $\Theta(s,1)$ take values in $[-\pi/2,\pi/2].$ In particular, one has $\Theta(0,0),\Theta(1,1)\in[-\pi/2,\pi/2].$ This proves part~(\Rn{2}).
\end{proof}

We are now ready to prove Proposition~\ref{prop_finiteTC}.

\begin{proof}[Proof of Proposition \ref{prop_finiteTC}]
For $a<b$, let
\[
\Omega_{a,b}=(a,b)\times(-1,1).
\]
By the Euler equation \eqref{Euler} and the slip boundary condition \eqref{SBC}, one has
\[
\partial_{x_2}P(\cdot,\pm1)=-\bv\cdot\nabla v_2(\cdot,\pm1)=0.
\]
Applying the divergence theorem, we obtain
\begin{align*}
\int_{\Omega_{a,b}}\Delta P\,dx
=-\int_{-1}^{1}\partial_{x_1}P(a,x_2)\,dx_2
+\int_{-1}^{1}\partial_{x_1}P(b,x_2)\,dx_2.
\end{align*}

Define $\widetilde{\msN}$ by
\begin{align}\label{eq_definition-tildeN}
\widetilde{\msN}=\msN\cup\{0,c\}\cup\msC(u)\cup\msC(u(a,\cdot))\cup\msC(u(b,\cdot)),    
\end{align}
where $\msC(u(a,\cdot))$ and $\msC(u(b,\cdot))$ denote the sets of critical values of the one-dimensional functions $u(a,\cdot)$ and $u(b,\cdot)$, respectively. By the Morse-Sard theorem, $\widetilde{\msN}$ is an $\mcL^1$-null subset of $\mbR$. For every $r\in u(\Omega_{a,b})\setminus\widetilde{\msN}$, condition (\Rn{3}) of Theorem \ref{thm_main} implies that $u^{-1}(r)$ must intersect the cross sections $\{a\}\times(-1,1)$ and $\{b\}\times(-1,1)$. Since $r\notin \msC(u(a,\cdot))\cup\msC(u(b,\cdot))$, at every such intersection point we have $\partial_{x_2}u\neq 0$. This means that $u^{-1}(r)$ intersects the cross sections transversely at those points. Consequently, the set $u^{-1}(r)\cap\overline{\Omega_{a,b}}\Subset\Omega_\infty$ consists of finitely many $C^2$ arcs, each of whose components intersects $\{a,b\}\times(-1,1)$ transversely.

For every $r\in u(\Omega_{a,b})\setminus\widetilde{\msN}$, let $\mcC_r(\Omega_{a,b})$ be the collection of connected components of $u^{-1}(r)\cap\overline{\Omega_{a,b}}$. For every $\Gamma\in\mcC_r(\Omega_{a,b})$, denote by $\omega_\Gamma$ the restriction of $\omega$ to $\Gamma$, which is a constant. Since $r$ is a regular value of $u$, there exists a simply connected neighborhood of $\Gamma$ on which $|\bv|>0$ and $\bv$ admits a lifting $\theta\in C^1$ with $\bv/|\bv|=(\cos\theta,\sin\theta)$. A direct computation shows that
\[
\nabla\theta=\frac{v_1\nabla v_2-v_2\nabla v_1}{|\bv|^2}.
\]
Denote by $\partial_s=\frac{\bv}{|\bv|}\cdot\nabla$ the tangential derivative along $\Gamma$, oriented by the flow. Now, applying the coarea formula yields
\begin{align*}
\int_{\Omega_{a,b}}\omega\bv\cdot\frac{v_1\nabla v_2-v_2\nabla v_1}{|\bv|^2}\,dx
=&\int_{\mbR\setminus\widetilde{\msN}}\int_{u^{-1}(r)\cap\Omega_{a,b}}\omega\partial_s\theta\,ds\,dr\\
=&\int_{\mbR\setminus\widetilde{\msN}}\sum_{\Gamma\in\mcC_r(\Omega_{a,b})}\omega_{\Gamma}\Delta_{\Gamma}\theta\,dr.
\end{align*}
Integrating \eqref{eq_key-identity} over $\Omega_{a,b}$, we arrive at
\begin{align}\label{eq_fundamental-formula}
\mcT(\bv,\Omega_{a,b})
=\int_{\mbR\setminus\widetilde{\msN}}\sum_{\Gamma\in\mcC_r(\Omega_{a,b})}\omega_{\Gamma}\Delta_{\Gamma}\theta\,dr
+\frac12\int_{-1}^{1}\partial_{x_1}P(a,x_2)\,dx_2
-\frac12\int_{-1}^{1}\partial_{x_1}P(b,x_2)\,dx_2.
\end{align}

For every $r\in u(\Omega_{a,b})\setminus\widetilde{\msN}$, denote by $N(r,\Omega_{a,b})$ the number of connected components of $u^{-1}(r)\cap\overline{\Omega_{a,b}}$, that is, the number of elements in $\mcC_r(\Omega_{a,b})$; for all other values of $r$, take $N(r,\Omega_{a,b})=0$. By Lemmas \ref{lemma_global-bound} and \ref{lemma_angle-bound},
\begin{align*}
\left|\sum_{\Gamma\in\mcC_r(\Omega_{a,b})}\omega_{\Gamma}\Delta_{\Gamma}\theta\right|
\le 2\pi\|\omega\|_{L^\infty(\Omega_\infty)}N(r,\Omega_{a,b}).    
\end{align*}
Since every $\Gamma\in\mcC_r(\Omega_{a,b})$ has exactly two endpoints at the vertical sides of $\partial\Omega_{a,b}$, we have
\begin{align*}
N(r,\Omega_{a,b})=\frac12\left( \#\{ u(a,x_2)=r \}+ \#\{ u(b,x_2)=r \} \right).
\end{align*}
Then by the area formula and Lemma \ref{lemma_global-bound},
\begin{align*}
\int_{\mbR}N(r,\Omega_{a,b})\,dr
=\frac12\int_{-1}^{1}(|\partial_{x_2}u(a,x_2)|+|\partial_{x_2}u(b,x_2)|)\,dx_2
\le 2\|\bv\|_{L^\infty(\Omega_\infty)}.
\end{align*}
Consequently,
\begin{align}\label{ineq_vorticity-angle}
\left|\int_{\mbR\setminus\widetilde{\msN}}\sum_{\Gamma\in\mcC_r(\Omega_{a,b})}\omega_{\Gamma}\Delta_{\Gamma}\theta\,dr\right|
\le 4\pi\|\bv\|_{L^\infty(\Omega_\infty)}\|\omega\|_{L^\infty(\Omega_\infty)}.
\end{align}

The horizontal component of the momentum equation reads
\[
v_1\partial_{x_1}v_1+v_2\partial_{x_2}v_1+\partial_{x_1}P=0.
\]
Integrating this equation with respect to $x_2\in[-1,1]$ and using the slip boundary condition \eqref{SBC} and the incompressibility condition, we obtain
\[
\frac{d}{d\,x_1}\int_{-1}^{1}(P(x_1,x_2)+v_1^2(x_1,x_2))\,dx_2=0.
\]
Hence, there exists a constant $C_0$ such that
\[
\int_{-1}^{1}P(x_1,x_2)\,dx_2=-\int_{-1}^{1}v_1^2(x_1,x_2)\,dx_2+C_0.
\]
This, together with Lemma \ref{lemma_global-bound}, implies that the $C^1$ function
\[
x_1\mapsto \int_{-1}^{1}P(x_1,x_2)\,dx_2
\]
is bounded. On the other hand, it is clear that a bounded $C^1$ function on a half-line has arbitrarily distant points where its derivative has arbitrarily small absolute value. Therefore, there exist a decreasing sequence $\{a_n\}$ with $a_n\to-\infty$ and an increasing sequence $\{b_n\}$ with $b_n\to+\infty$ such that
\begin{align}\label{eq_limit-dx1P}
\left|\int_{-1}^{1}\partial_{x_1}P(a_n,x_2)\,dx_2\right|
+\left|\int_{-1}^{1}\partial_{x_1}P(b_n,x_2)\,dx_2\right|\longrightarrow0.
\end{align}

Finally, choose $a=a_n$ and $b=b_n$ in \eqref{eq_fundamental-formula}, and redefine $\widetilde{\msN}$ in \eqref{eq_definition-tildeN} by
\[
\widetilde{\msN}=\msN\cup\{0,c\}\cup\msC(u)\cup \bigcup_{n\ge1}\big(\msC(u(a_n,\cdot))\cup\msC(u(b_n,\cdot))\big),
\]
which is still $\mcL^1$-null. Passing to the limit as $n\to\infty$ and using \eqref{ineq_vorticity-angle} and \eqref{eq_limit-dx1P}, we obtain \eqref{ineq_TCbound}. This completes the proof.
\end{proof}



\medskip

\section{Asymptotic behavior}\label{sec_asymptotic}

In this section, equipped with Lemma \ref{lemma_global-bound} and Proposition \ref{prop_finiteTC}, we study the asymptotic behavior of the flow in the far field.

Throughout this section, let $\bv=\nabla^\perp u\in C^1(\overline{\Omega_\infty})$ and $P\in C^1(\overline{\Omega_\infty})$ solve \eqref{Euler}--\eqref{SBC} such that assumption (\Rn{3}) of Theorem~\ref{thm_main} is fulfilled.

First, since $\bv$ is bounded by Lemma \ref{lemma_global-bound} and the total curvature is finite by Proposition \ref{prop_finiteTC}, the flow satisfies the assumption in Theorem~1.1 of \cite{GLXX_arXiv_2026}. Consequently, the pressure $P$ converges in the far field.
\begin{lemma}[cf.\ \cite{GLXX_arXiv_2026}]\label{lemma_limitP}
The pressure $P$ is bounded, and there exist two numbers $P_+$ and $P_-$ such that
\begin{align*}
\lim_{x_1\to\pm\infty}P(x_1,x_2)=P_\pm \quad \hbox{ uniformly in } x_2\in[-1,1].
\end{align*}    
\end{lemma}

The convergence of $P$ in the far field, together with the global bounds in Lemma \ref{lemma_global-bound}, implies the asymptotic vanishing of the vertical velocity component.

\begin{lemma}\label{lemma_v2converge0}
It holds that
\begin{align*}
\lim_{|x_1|\to\infty}v_2(x_1,x_2)=0,\quad \hbox{ uniformly in } x_2\in[-1,1].
\end{align*}
\end{lemma}
\begin{proof}
We prove the lemma by contradiction. Assume that there exist $\varepsilon>0$, a sequence $\{a_n\}$ with $|a_n|\to+\infty$, and numbers $b_n\in[-1,1]$ such that
\begin{equation*}
|v_2(a_n,b_n)|\ge \varepsilon.
\end{equation*}   
Passing to subsequences if necessary, we may assume without loss of generality that
\[
a_n\to+\infty,\quad b_n\to b_*\in[-1,1].
\]
Define
\[
\bv^{(n)}(x_1,x_2)=\bv(a_n+x_1,x_2),
\quad
P^{(n)}(x_1,x_2)=P(a_n+x_1,x_2)-P_+,
\]
where $P_+$ is the limit of $P$ at $x_1=+\infty$, as in Lemma \ref{lemma_limitP}.

Fix a number $p\in(2,\infty)$ and let $\alpha=1-\frac2p$. Fix temporarily $k\in\mbN_+$ and let $\Omega_k=(-k,k)\times(-1,1)$.
It follows from Lemma \ref{lemma_limitP} that
\begin{align}\label{eq_Pn-converges-0}
P^{(n)}\longrightarrow 0 \quad \hbox{in } C(\overline{\Omega_k}).    
\end{align}
By Lemma \ref{lemma_global-bound}, $\bv^{(n)}$ is bounded uniformly in $n$ in $W^{1,p}(\Omega_k)$; hence, it is uniformly bounded in $C^{0,\alpha}(\overline{\Omega_k})$ by Morrey's inequality.
Consequently, up to a subsequence, there exists a vector field $\bw=(w_1,w_2)\in W^{1,p}(\Omega_k)$ such that
\[
\bv^{(n)}\rightharpoonup \bw \quad \hbox{weakly in } W^{1,p}(\Omega_k),
\]
and
\begin{align}\label{eq_vn-converges-w}
\bv^{(n)} \longrightarrow \bw \quad \hbox{in } C(\overline{\Omega_k}).    
\end{align}
In particular, one has
\begin{equation}\label{ineq_w2-lower-bound}
|w_2(0,b_*)|
=\lim_{n\to\infty}|v_2(a_n,b_*)|=\lim_{n\to\infty}|v_2(a_n,b_n)|
\ge\varepsilon.
\end{equation}
Since $k\in\mbN_+$ is arbitrary, a diagonal argument yields a vector field $\bw=(w_1,w_2)\in W^{1,p}_{\rm loc}(\overline{\Omega_\infty})$ and a subsequence (still denoted by the same notation) such that \eqref{eq_Pn-converges-0} and \eqref{eq_vn-converges-w} hold for every $k\in\mbN_+$.

On the other hand, since $(\bv^{(n)},P^{(n)})$ solves \eqref{Euler}-\eqref{SBC}, by \eqref{eq_Pn-converges-0} and \eqref{eq_vn-converges-w}, we can pass to the limit in the weak form of \eqref{Euler} to obtain
\begin{eqnarray*}
\left\{\begin{aligned}
&\Div\,(\bw\otimes\bw)=0, & \hbox{in } \Omega_\infty, \\
&\Div\,\bw=0, & \hbox{in } \Omega_\infty, \\
&w_2=0, & \hbox{on } \partial\Omega_\infty.
\end{aligned}\right.
\end{eqnarray*}
However, any solution to this system in the class $H^1_{\rm loc}(\overline{\Omega_\infty})\cap C(\overline{\Omega_\infty})$ must satisfy $w_2=0$ in $\Omega_\infty$ (see, e.g., \cite{GXX_CMP_2026}), which contradicts \eqref{ineq_w2-lower-bound}. This completes the proof of the lemma.
\end{proof}

\begin{lemma}\label{lemma_sequence-u-converge}
For any sequences $\{a_n\}$ with $a_n\to-\infty$ and $\{b_n\}$ with $b_n\to+\infty$, there exist subsequences (still denoted by $\{a_n\}$ and $\{b_n\}$) and functions $U_-,U_+\in C^{1}([-1,1])$ such that
\begin{align*}
u(a_n+\cdot,\cdot)\longrightarrow U_-,\quad u(b_n+\cdot,\cdot)\longrightarrow U_+ \quad \hbox{in } C^1([-k,k]\times[-1,1])
\end{align*}
for every $k\in\mbN_+$.
\end{lemma}
\begin{proof}
For every fixed $k\in \mbN_+$ and $p>2$, by Lemma~\ref{lemma_global-bound}, the translates $u(a_n+\cdot,\cdot)$ and $u(b_n+\cdot,\cdot)$ are uniformly bounded in $W^{2,p}((-k,k)\times(-1,1))$. The compact embedding
\[
W^{2,p}((-k,k)\times(-1,1))\Subset C^1([-k,k]\times[-1,1])
\]
and a diagonal extraction give the desired subsequences converging in $C^1([-k,k]\times[-1,1])$ for every $k\in\mbN_+$. Lemma~\ref{lemma_v2converge0} then identifies the limits $U_\pm$ as functions of $x_2$ alone. This completes the proof.
\end{proof}

Note that the limiting profiles $U_\pm$ are selected limits and are not assumed to be unique.


\medskip

\section{Vanishing of the total curvature}\label{sec_vanishTC}

This final section is devoted to the proof of Theorem \ref{thm_main}.

For simplicity of notation, define
\[
F(x_1)=\frac12\int_{-1}^{1}\partial_{x_1}P(x_1,x_2)\,dx_2.
\]
By \eqref{eq_limit-dx1P} and Lemma \ref{lemma_sequence-u-converge}, there exist a decreasing sequence $\{a_n\}$ with $a_n\to-\infty$ and an increasing sequence $\{b_n\}$ with $b_n\to+\infty$ such that
\begin{align}\label{eq_limitFn}
|F(a_n)|+|F(b_n)|\longrightarrow0,
\end{align}
and
\begin{align}\label{eq_sequence-u-converge}
u(a_n+\cdot,\cdot)\longrightarrow U_-,\quad u(b_n+\cdot,\cdot)\longrightarrow U_+ \quad \hbox{in } C^1([-k,k]\times[-1,1])
\end{align}
for every $k\in\mbN_+$. Henceforth, we fix these sequences and denote
\[
\Omega_n=(a_n,b_n)\times(-1,1).
\]
Redefine $\widetilde{\msN}$ in \eqref{eq_definition-tildeN} by
\[
\widetilde{\msN}=\msN\cup\{0,c\}\cup\msC(u)\cup\msC(U_+)\cup\msC(U_-)\cup\bigcup_{n\ge1}\big(\msC(u(a_n,\cdot))\cup\msC(u(b_n,\cdot))\big),
\]
which is an $\mcL^1$-null subset of $\mbR$ by Morse-Sard theorem. Let $\widetilde{\msR}=u(\Omega_\infty)\setminus\widetilde{\msN}$. Then the following functions of $r$ are defined for $r\in\widetilde{\msR}$ and take the value zero for $r\in\mbR\setminus\widetilde{\msR}$. Define
\begin{align*}
T_n(r)=\sum_{\Gamma\in\mcC_r(\Omega_n)}|\Delta_{\Gamma}\theta|,
\quad
N_n(r)=\#\mcC_r(\Omega_n),
\end{align*}
where $\mcC_r(\Omega_n)$ denotes the collection of connected components of $u^{-1}(r)\cap \overline{\Omega_{n}}$. Remember that
\[
N_n(r)=\frac12\left( \#\{ u(a_n,x_2)=r \}+ \#\{ u(b_n,x_2)=r \} \right).
\]
We therefore introduce
\[
N_\infty(r)=\frac12\left( \#(U_-)^{-1}(r)+ \#(U_+)^{-1}(r)\right),
\]
and expect $N_n$ to converge to $N_\infty$ in some sense.

It now follows from \eqref{eq_fundamental-formula} that
\begin{align}\label{ineq_TCn-bound}
\mcT(\bv,\Omega_{n})
\le \|\omega\|_{L^\infty(\Omega_\infty)}\int_{\mbR}T_n(r)\,dr
+|F(a_n)|+|F(b_n)|.
\end{align}
By Lemma \ref{lemma_angle-bound},
\begin{align}\label{ineq_2nd-condition}
T_n(r)\le 2\pi N_n(r).
\end{align}
In view of \eqref{eq_limitFn} and \eqref{ineq_TCn-bound}, the proof of Theorem \ref{thm_main} reduces to showing that
\[
\int_{\mbR}T_n(r)\,dr\longrightarrow0.
\]
For this, we shall use the following elementary lemma, which is a substitute for dominated convergence theorem.
\begin{lemma}\label{lemma_gDCT}
Let $\{f_n\}$ and $\{g_n\}$ be two sequences of nonnegative functions in $L^1(\mbR)$. Assume that 
\begin{enumerate}[label=\textup{(\roman*)}]
\item $f_n\longrightarrow0$ a.e. in $\mbR$;

\item for every $n$, $f_n\le g_n$ a.e. in $\mbR$;

\item there exists a function $g_\infty\in L^1(\mbR)$ such that  $g_n\longrightarrow g_\infty$ in $L^1(\mbR)$.
\end{enumerate}
Then it holds that
\[
\lim_{n\to\infty}\int_{\mbR}f_n\,dr=0.
\]
\end{lemma}
\begin{proof}
By assumption (\Rn{2}), one has
\[
f_n\le \min\{f_n, g_\infty \}+|g_n-g_\infty|, \quad \hbox{a.e. in } \mbR. 
\]
Note that $\min\{f_n, g_\infty \}$ is bounded by $g_\infty\in L^1(\mbR)$, and by assumption (\Rn{1}), converges to zero a.e. in $\mbR$. So by the dominated convergence theorem,
\[
\lim_{n\to\infty}\int_{\mbR}\min\{f_n, g_\infty \}\,dr=0.
\]
This, together with assumption (\Rn{3}), gives the desired result.
\end{proof}

We shall apply Lemma \ref{lemma_gDCT} by choosing $f_n=T_n(r)$ and $g_n=2\pi N_n(r)$. The second condition in Lemma \ref{lemma_gDCT} has been verified by \eqref{ineq_2nd-condition}, so it remains to verify the first and the third conditions. The following elementary lemma will be needed for both purposes.

\begin{lemma}\label{lemma_rootnumber}
Let $h_n\to h$ in $C^1([-1,1])$, and let $r$ be a regular value of $h$ distinct from $h(-1)$ and $h(1)$. Then $h^{-1}(r)$ is finite and
\[
 \#h_n^{-1}(r)=\#h^{-1}(r)
\]
for all sufficiently large $n$.
Moreover, there exists $c_0=c_0(r)>0$ such that for all sufficiently large $n$, every $t\in h_n^{-1}(r)$ satisfies $|h_n'(t)|\ge c_0$.
\end{lemma}
\begin{proof}
Since $r\neq h(\pm1)$, we have $h^{-1}(r)\subset(-1,1)$. If $a\in h^{-1}(r)$, then $h'(a)\neq0$; hence $a$ is an isolated point of $h^{-1}(r)$. By compactness, $h^{-1}(r)$ is finite. Now write
\[
h^{-1}(r)=\{a_1<\cdots<a_m\}.
\]
Define
\[
m_0:=\min_{1\le j\le m}|h'(a_j)|>0.
\]
Choose pairwise disjoint intervals
\[
I_j=(a_j-\delta_j,a_j+\delta_j)\subset(-1,1)
\]
such that on each $\overline{I_j}$,
\[
|h'|\ge \frac{m_0}{2}
\]
and $h'$ has constant sign. It is clear that
\[
\rho:=\min_{x\in J}|h(x)-r|>0, \quad \hbox{with } J:=[-1,1]\setminus\cup_{j=1}^m I_j,
\]
and
\[
\eta:=\min_{1\le j\le m}
\min\{|h(a_j-\delta_j)-r|,\ |h(a_j+\delta_j)-r|\}>0.
\]
Choose $n_0$ large enough such that for all $n\ge n_0$,
\[
\|h_n-h\|_{L^\infty(-1,1)}<\min\left\{\frac{\rho}{2},\frac{\eta}{2}\right\},
\quad
\|h_n'-h'\|_{L^\infty(-1,1)}<\frac{m_0}{4}.
\]
Then $h_n\neq r$ on $J$. On each $I_j$, the endpoint values of $h_n-r$ have the same signs as those of $h-r$, hence opposite signs. Moreover, on $I_j$,
\[
|h_n'|\ge\frac{m_0}{4}>0.
\]
Hence $h_n$ is strictly monotone on $I_j$, and $h_n-r$ has exactly one zero in $I_j$. Thus, for all $n\ge n_0$,
\[
\#h_n^{-1}(r)=m=\#h^{-1}(r),
\]
and every $t\in h_n^{-1}(r)$ satisfies
\[
|h_n'(t)|\ge \frac{m_0}{4}.
\]
This completes the proof.
\end{proof}

\begin{lemma}\label{lemma_Tn-pointwise-converge}
It holds that
\begin{align*}
T_n(r)\longrightarrow 0,\quad \hbox{for a.e. } r\in\mbR. 
\end{align*}
\end{lemma}
\begin{proof}
Fix any $r\in\widetilde{\msR}$.

First, we show that $|v_1|$ has a positive lower bound on $u^{-1}(r)$ in the far field. By \eqref{eq_sequence-u-converge} and Lemma \ref{lemma_rootnumber}, there exists $c_0(r)>0$ such that, for sufficiently large $n$, on the intersection of $u^{-1}(r)$ with the cross sections $\{a_n\}\times(-1,1)$ and $\{b_n\}\times(-1,1)$,
\begin{align}\label{ineq_v1-lower-bound}
|v_1|=|\partial_{x_2}u|\ge c_0(r).    
\end{align}
Let $\Gamma$ be any connected component of $u^{-1}(r)$. By Bernoulli's principle, there exists a constant $H_\Gamma$ such that
\[
|\bv|^2=2(H_\Gamma-P) \quad \hbox{on } \Gamma.
\]
By Lemma \ref{lemma_limitP} and \eqref{ineq_v1-lower-bound}, we have
\[
2(H_\Gamma-P_+)\ge c_0^2(r), \quad 2(H_\Gamma-P_-)\ge c_0^2(r).
\]
Then for $x\in\Gamma$ with $|x_1|$ sufficiently large, by Lemmas \ref{lemma_limitP} and \ref{lemma_v2converge0}, 
\begin{align}\label{ineq_v1-true-lower-bound}
v_1^2=2(H_\Gamma-P)-v_2^2\ge\frac12 c_0^2(r).
\end{align}
This lower bound, while depending on $r$, can be made independent of $\Gamma$, since $u^{-1}(r)$ has finitely many components by assumption (\Rn{3}) of Theorem~\ref{thm_main}.

Combining \eqref{ineq_v1-true-lower-bound} and Lemma \ref{lemma_v2converge0}, we infer that there exists $L(r)>0$ such that, for all $a,b$ with $a<-L(r)<L(r)<b$, the intersection of every component of $u^{-1}(r)$ with $[a,b]\times[-1,1]$ is exactly one compact arc joining the two cross sections.

Now, for $n$ sufficiently large, every connected component $\Gamma$ of $u^{-1}(r)\cap\Omega_n$ has two endpoints $(a_n,c_n)$ and $(b_n,d_n)$ for some $c_n,d_n\in[0,1]$. Applying Lemma \ref{lemma_angle-bound} gives
\begin{align*}
|\Delta_\Gamma\theta|=&\left|\arctan\left( \frac{v_2(a_n,c_n)}{v_1(a_n,c_n)}\right)- \arctan\left( \frac{v_2(b_n,d_n)}{v_1(b_n,d_n)}\right)\right|.
\end{align*}
By \eqref{ineq_v1-true-lower-bound} and Lemma \ref{lemma_v2converge0} again, it is evident that the right-hand side converges to zero. Since the number of connected components of $u^{-1}(r)$ is finite, we see that $T_n(r)\to0$. This completes the proof.
\end{proof}

\begin{lemma}\label{lemma_Nn-converge}
It holds that
\begin{align*}
N_n\longrightarrow N_\infty, \quad \hbox{in } L^1(\mbR). 
\end{align*}
\end{lemma}
\begin{proof}
First, by the area formula,
\begin{align*}
\int_{\mbR}N_n(r)\,dr
=\frac12\int_{-1}^{1}(|\partial_{x_2}u(a_n,x_2)|+|\partial_{x_2}u(b_n,x_2)|)\,dx_2.
\end{align*}
Then, by \eqref{eq_sequence-u-converge} and the area formula again,
\begin{align*}
\lim_{n\to\infty}\int_{\mbR}N_n(r)\,dr
=&\frac12\int_{-1}^{1}(|U_{-}^{'}(x_2)|+|U_{+}^{'}(x_2)|)\,dx_2\\
=&\frac12\int_{\mbR}\left( \#(U_-)^{-1}(r)+ \#(U_+)^{-1}(r)\right) \,dr\\
=&\int_{\mbR} N_\infty(r) \,dr.
\end{align*}
Second, by \eqref{eq_sequence-u-converge} and Lemma \ref{lemma_rootnumber},
\begin{align*}
N_n(r)\longrightarrow N_\infty(r),\quad \hbox{ for a.e. } r\in\mbR.
\end{align*}
Finally, since both $N_n$ and $N_\infty$ are nonnegative, applying Scheff\'e's lemma finishes the proof of the lemma.
\end{proof}

We conclude with the proof of Theorem \ref{thm_main}.

\begin{proof}[Proof of Theorem \ref{thm_main}]
We only need to show that (\Rn{3}) implies (\Rn{1}). By \eqref{ineq_2nd-condition} and Lemmas \ref{lemma_gDCT}, \ref{lemma_Tn-pointwise-converge}, and \ref{lemma_Nn-converge},
\[
\lim_{n\to\infty}\int_{\mbR}T_n(r)\,dr=0.
\]
This, together with \eqref{eq_limitFn} and \eqref{ineq_TCn-bound}, implies that $\mcT(\bv,\Omega_\infty)=0$. By the slip boundary condition, it then readily follows that $\bv$ is a shear flow of the form $\bv=(v_1(x_2),0)$ (see, e.g., \cite{GXX_CMP_2026}). This completes the proof of Theorem \ref{thm_main}.
\end{proof}


\medskip

\noindent {\bf Acknowledgments.}
The research of Gui is supported by the NSFC Key Program (Grant No.~12531010), University of Macau research grants CPG2024-00016-FST, CPG2025-00032-FST, SRG2023-00011-FST, MYRG-GRG2023-00139-FST-UMDF, the UMDF Professorial Fellowship of Mathematics, Macao SAR FDCT 0003/2023/RIA1, and Macao SAR FDCT 0024/2023/RIB1. The research of Yang is partially supported by National Key R$\&$D Program of China 2022YFA1006800, NSFC No.~12531010, FDCT No.~0070/2024/RIA1, Multi-Year Research Grant No.~MYRG-GRG2024-00082-FST-UMDF, MYRGGRG202500051-FST and UMDF No.~TISF/2025/006/FST. The authors would like to express their sincere gratitude to Chunjing Xie for his generous and valuable discussions on this work.

\medskip

\noindent {\bf Disclosure of AI use.} 
The authors used the 5.6 Sol and 6.1 Sol models in ChatGPT Plus (OpenAI) to discuss selected technical aspects of this work. The authors also used 6.1 Sol for proofreading shortly before this paper was posted online. The research problem, main ideas, and overall proof strategy were developed by the authors, building on their previous work on steady Euler flows. The manuscript was written entirely by the authors and was not generated by a large language model. The authors take full responsibility for all mathematical arguments and the final text.


\end{document}